\documentclass[11pt]{amsart}

\usepackage{amssymb} % mathematical fonts (in amsfonts)
\usepackage[mathscr]{eucal} % load euler fonts (in amsfonts)
\usepackage[all]{xy}
\usepackage{amsmath} %this is needed for \rxar (\xrightarrow)
\usepackage{epsfig}
\usepackage{amscd}

\newtheorem{TEO}{Theorem}[section]

\newtheorem{LEM}[TEO]{Lemma}

\newtheorem{COR}[TEO]{Corollary}

\newtheorem{CLA}[TEO]{Claim}

\newtheorem{theoremalpha}{Theorem}

\newtheorem*{introcorollary}{Corollary}

\theoremstyle{definition}
\newtheorem{remark}[TEO]{Remark}

\newcommand\Oh{{\mathcal O}}

\def\OO{{\mathcal O}}

\newcommand\dual{\mathrel{\raise3pt\hbox{$\underline{\mathrm{\thinspace d
\thinspace}}$}}}

\newcommand\proj{\mathbb P}

\newcommand\coker{\operatorname{coker}}

\newcommand\Pico{\operatorname{Pic^0}}

\begin{document}

\footnote{
Partially supported by PRIN 2007 MIUR: "Spazi dei moduli e teoria di Lie" and PRIN 2006 of MIUR "Geometry on algebraic varieties".

 }

\title{Hyperplane sections of abelian surfaces }

\author[E. Colombo]{Elisabetta Colombo}
\address{Dipartimento di Matematica,
Universit\`a di Milano, via Saldini 50,
     I-20133, Milano, Italy } \email{{\tt
elisabetta.colombo@unimi.it}}

\author[P. Frediani]{Paola Frediani}
\address{ Dipartimento di Matematica, Universit\`a di Pavia,
via Ferrata 1, I-27100 Pavia, Italy } \email{{\tt
paola.frediani@unipv.it}}

\author[G. Pareschi]{Giuseppe Pareschi}
\address{Dipartimento di Matematica, Universit\`a di Roma, Tor
Vergata, V.le della Ricerca Scientifica, I-00133 Roma, Italy}
\email{{\tt pareschi@mat.uniroma2.it}}

%\date{\today}
\maketitle

\setlength{\parskip}{.1 in}

\begin{abstract}
By a theorem of Wahl, the canonically embedded curves which are
hyperplane section of K3 surfaces are distinguished by the
non-surjectivity of their Wahl map. In this paper we address the
problem of distinguishing  hyperplane sections of abelian surfaces.
The somewhat surprising result is that the Wahl map of such curves
is (tendentially) surjective, but their \emph{second} Wahl map has
corank at least 2 (in fact a more precise result is proved).
\end{abstract}

\section{Introduction}
Which canonically embedded curves $C\subset \proj^{g-1}$ are
hyperplane sections of K3 surfaces? An answer to this question is
provided by Wahl's theorem (\cite{wahl1},\cite{wahl3} see also the
proof by Beauville-M\'erindol in \cite{bm}), which asserts that, if
$C$ is a curve sitting on a $K3$ surface, then its Wahl map, or
gaussian map,
$$\gamma^1_C:\wedge^2H^0(K_C)\rightarrow H^0(K_C^3)$$
is not surjective. The significance of this criterion is better
appreciated     if one compares it with two other results:

 \noindent - Ciliberto-Harris-Miranda's
theorem (\cite{chm}, see also Voisin's proof in \cite{voi}), stating
that the Wahl map of the generic curve of genus $g$ is surjective as
soon as this is numerically possible, i.e. for $g\ge 10$, with the
exception of $g=11$. For $g<10$ and $g=11$ it is  known that the
generic curve lies on a K3 surface (\cite{mm}).

\noindent
 - Lazarsfeld theorem (\cite{laz}, see also \cite{pa} for a  different proof), asserting (loosely speaking) that
general hyperplane sections of general polarized $K3$ surfaces
satisfy the Brill-Noether-Petri condition. Hence curves lying on K3
surfaces cannot be distinguished by special Brill-Noether-theoretic
properties. In fact it has been independently conjectured by various
authors (Mukai, Voisin, see \cite{voi}4.13(b), and Wahl, see
\cite{wahl4}\S0) that for Brill-Noether-Petri-general curves the
non-surjectivity of the Wahl map should completely characterize
curves contained in   K3 surfaces.

The aim of this note is to address the same questions for the other
class of canonically embedded hyperplane sections of smooth
surfaces: hyperplane sections of abelian surfaces. This case is
somewhat more subtle since here the canonical embedding
$C\hookrightarrow \proj^{g-3}$ is not complete, as it is obtained by
the complete canonical embedding in $ \proj^{g-1}=\proj\,
H^1(\Oh_C)$ after projection from the line $\proj H^1( \Oh_X)\subset
\proj H^1(\Oh_C)$ (where $X$ is the abelian surface). This is the
first-order infinitesimal counterpart of the obvious intrinsic
restriction satisfied by a curve sitting on an abelian surface,
namely that its Jacobian is non-simple, as it contains the abelian
surface $\Pico X$. However, given a canonically embedded curve
$C\subset \proj^{g-3}$, obtained by projection from a line $L$ in
$\proj^{g-1}$, it is not obvious how to recognize that its jacobian
contains a two-dimensional abelian surface $Y$, so that $L$ is the
projectivized tangent space to $Y$. This is one of the features of
our main result.

In order to set the stage, let us recall that the Wahl map
$\gamma^1$ belongs to a hierarchy of maps, called higher gaussian
maps. The second gaussian map, or second Wahl map, is a linear map
$$\gamma^2_C: I_2(C)\rightarrow H^0(K_C^4)$$
where $I_2(C)$ is the kernel of the natural map $S^2H^0(K_C)
\rightarrow H^0(K_C^2)$ i.e. (if $C$ is non-hyperelliptic) the
vector space of quadrics of $\proj^{g-1}=\proj\; H^1(\Oh_C)$
containing $C$. The map $\gamma^2_C$ has an independent interest,
which was the original motivation of the first two authors for
studying this sort of questions.  In fact there is a relation,
analyzed  in \cite{cf2}, between the second gaussian map and the
curvature of the moduli space ${M_g}$ of curves of genus $g$,
endowed with the Siegel metric induced by the period map $j:{
M_g}\rightarrow {A_g}$. To be precise, in \cite{cf2}  the
holomorphic sectional curvature of ${M_g}$ along the
 a Schiffer variation $\xi_P$, for $P$ a point  on the curve $C$,
 was computed in terms of the holomorphic sectional curvature of ${A_g}$ and
the
 map $\gamma^2_C$. This was accomplished using the formula for the
second fundamental form associated to the period map given in
\cite{cpt}.

Going back to our problem, we introduce the following notation.
Given a  subspace $W\subset H^0(K_C)$, we will denote
$$S^2W\cdot H^0(K_C^2)$$ the image of $S^2W\otimes H^0(K_C^2)$ in $H^0(K_C^4)$ via
the natural multiplication map. If $W$ is $2$-dimensional and base-point-free, the base-point-free
pencil trick implies that $S^2W\cdot H^0(K_C^2)$ has codimension 2
in $H^0(K_C^4)$. If $C$ is embedded in abelian surface, then
$H^0(\Omega^1_X)$ is naturally a (base-point-free) $2$-dimensional subspace of
$H^0(K_C)$.
 Our main
result is
\begin{theoremalpha}\label{main}
Let $C$ be a curve contained in abelian surface $X$. Then the image
of the second gaussian map $\gamma^2_C$ is contained in
$S^2H^0(\Omega^1_X)\cdot H^0(K_C^2)$ (notation as above). Therefore
the corank of $\gamma^2_C$ is at least 2.

Moreover, if the second gaussian map of the surface $X$ \emph{(see
\S2)} is surjective, then the image of the map $\gamma^2_C$ concides
with $S^2H^0(\Omega^1_X)\cdot H^0(K_C^2).$
\end{theoremalpha}
The above Theorem can be stated, perhaps more suggestively, as
follows. Given a subspace $V\subset H^1(\OO_C)$, let $\overline
V\subset H^0(K_C)$ be its conjugate.
\begin{introcorollary}\label{version} Let $C\subset \proj^{g-3}$ be a canonically embedded
curve of genus $g$, obtained from the complete canonical embedding
$C\subset \proj H^1(\OO_C)=\proj^{g-1}$ by projection from a line
$\proj V\subset \proj H^1(\OO_C)$,  $\dim V=2$. If $C$ is a
hyperplane section of an abelian surface $X\subset \proj^{g-2}$ then
$$Im(\gamma^2_C)\subseteq S^2\overline{V}\cdot H^0(K_C^2)$$
\end{introcorollary}
A few comments are in order. In the first place,
 the \emph{first} gaussian map of a curve $C$ sitting
on an abelian surface $X$ is "tendentially surjective". This is the
content of another result proved in this note (we refer to Theorem
\ref{firstgaussian} for the precise statement):
\begin{theoremalpha}\label{loose}  Assume that the first
gaussian map of the line bundle $\OO_X(C)$ on the surface $X$ is
surjective and that the multiplication map
$$\gamma^0_{X,C}:H^0(X,\OO_X(C))\otimes H^0(C,K_C)\rightarrow H^0(K_C^2)$$
is surjective
\emph{(for example, both conditions hold if $\OO_X(C)$ is at
least a $5$-th power of a ample line bundle on $X$, see
\cite{par})}. Then the first gaussian map of $C$ is surjective.
\end{theoremalpha}

Secondly, one expects a Ciliberto-Harris-Miranda's theorem for
second gaussian maps, namely that, for the generic curve of genus
$g\ge 18$, the second gaussian map $\gamma^2$ is surjective. In
\cite{cf1} the first two authors exhibited infinitely many genera
where this happens, by producing examples of curves lying on the
product of two curves with surjective second Gaussian map whose
second Gaussian map is surjective. Other examples were given in
\cite{bafo}. Both classes of examples generalize constructions
given by Wahl for the first Gaussian map in \cite{wahl3},
\cite{wahl1}. Moreover the first two authors (\cite{cf3}) have
proved the surjectivity of $\gamma^2$ for the generic curve of
high genus (for $g> 152$).

Finally, concerning the Brill-Noether theory of curves on abelian
surfaces, M. Paris (\cite{paris}) has obtained the following almost
complete extension of Lazarsfeld's result: let $X$ be an abelian
surface such that its N\'eron-Severi group $NS(X)$ is cyclic,
spanned by $c_1(L)$, with $L$ an ample line bundle on $X$. Then all
line bundles of degree $d\ne g(C)-1$ on a general curve $C\in |L|$
satisfy the Petri condition. Hence, at least in degree different
from $g-1$, there are no Brill-Noether-theoretic ways to distinguish
curves sitting on abelian surfaces. It is reasonable to conjecture
that -- under the hypothesis of sufficient Brill-Noether generality
-- the conclusion of the Corollary  should characterize hyperplane
section of abelian surfaces.

The proofs of both Theorems \ref{main} and \ref{loose} are based on
cohomological computations concerning the extension classes of the
cotangent sequence
$$0\rightarrow K_C^{-1}\rightarrow {\Omega^1_X}_{|C}\rightarrow
K_C\rightarrow 0$$ and of its "symmetric square" $$0\rightarrow
{\Omega^1_X}_{|C}\otimes K_C^{-1}\rightarrow
S^2{\Omega^1_X}_{|C}\rightarrow K_C^2\rightarrow 0$$ Our approach
has its roots in Beauville-M\'erindol's paper (\cite{bm}).

Finally, we remark that, beyond Wahl's theorem, the geometric
significance of Wahl's map is now reasonably well understood, thanks
to the work of many authors (notably Voisin, \cite{voi}). On the
other hand, Theorem \ref{main} above, as well as the works
\cite{cf1},\cite{cf2}, seem to indicate that also the second
gaussian map encodes some interesting geometry, which is at present
much less understood.

The third author thanks Andrea Susa for his invaluable help.

\section{Preliminaries on gaussian maps}

\subsection{Classical gaussian maps (of any order)}
 Let $Y$ be a smooth complex projective variety and let $\Delta_Y\subset
Y\times Y$ be the diagonal. Let $L$ and $M$ be line bundles on $Y$.
For a non-negative integer $k$, the \emph{k-th Gaussian map}
associated to these data is the restriction map to the diagonal
\begin{equation}\label{gaussian1}\gamma^k_{L,M}:H^0(Y\times Y,I^k_{\Delta_Y}\otimes
L\boxtimes M )\rightarrow
H^0(Y,{I^k_{\Delta_Y}}_{|\Delta_Y}\otimes L\otimes M)\cong
H^0(Y,S^k\Omega_Y^1\otimes L\otimes M).
\end{equation}
Usually  \emph{first} gaussian maps are simply referred to as
 \emph{gaussian maps}. The exact sequence
\begin{equation}
\label{Ik} 0 \rightarrow I^{k+1}_{\Delta_Y} \rightarrow
I^k_{\Delta_Y} \rightarrow S^k\Omega^1_Y \rightarrow 0,
\end{equation}(where $S^k\Omega^1_Y$ is identified to its image via the diagonal map), twisted by $L\boxtimes M$, shows that the domain of the $k$-th
gaussian map is the kernel of the previous one:
$$\gamma^k_{L,M}:
ker \gamma^{k-1}_{L,M}\rightarrow H^0(S^k\Omega_Y^1\otimes L\otimes
M).$$ In our applications, we will exclusively deal with gaussian
maps of order one and two, assuming also that the two line bundles
$L$ and $M$ coincide. For the reader's convenience, we spell out
these maps. The map $\gamma^0_L$ is the multiplication map of global
sections
\begin{equation}\label{symmetric}H^0(X,L)\otimes
H^0(X,L)\rightarrow H^0(X,L^2)\end{equation}
 which obviously
vanishes identically on $\wedge^2 H^0(L)$.
 Consequently, $H^0(Y
\times Y, I_{\Delta_Y}\otimes L\boxtimes L)$ decomposes as
$\wedge^2 H^0(L)\oplus I_2(L)$, where $I_2(L)$ is the kernel of
$S^2H^0(X,L)\rightarrow H^0(X,L^2)$. Since $\gamma^1_L$ vanishes
on symmetric tensors, one writes
\begin{equation}\label{gamma1L}\gamma^1_L:\wedge^2H^0(L)\rightarrow H^0(\Omega^1_Y\otimes
L^2).\end{equation}
 Again, $H^0(Y\times Y, I^2_{\Delta_Y}\otimes L\boxtimes L)$
 decomposes as  the sum of $I_2(L)$ and the kernel of (\ref{gamma1L}). Since $\gamma_L^2$ vanishes identically on skew-symmetric
 tensors, one usually writes
 \begin{equation}\label{skew}\gamma^2_L:I_2(L)\rightarrow H^0(S^2\Omega_Y^1\otimes L^2)
 \end{equation}
 (In general, gaussian maps of even (resp.odd) order vanish
 identically on skew-symmetric (resp. symmetric) tensors). The primary object of this paper will be the
 \emph{first and second order Wahl maps}, which
 are the first and second gaussian maps of the canonical line bundle
 $K_C$ on a curve $C$:
 $$\gamma^1_C:\wedge^2H^0(K_C)\rightarrow
H^0(K_C^3)$$
 $$\gamma^2_C:I_2(K_C)\rightarrow H^0(K_C^4)$$

 \subsection{Curve-surface gaussian maps}
In general, given a variety $Y$, endowed with a divisor $Z$ on $Y$,
it is
 useful to consider a variant of gaussian maps, which lies somewhat
 in between the gaussian maps on $Y$ and the ones on $Z$. When $Y$ is a surface and $C$ is a curve on it,
 for easy reference we will
 sometimes call them \emph{curve-surface gaussian maps}. These maps already appear in \cite{voi}. They are simply defined
 as follows: let $L$ be a line bundle on $Y$, and let $M_Z$ be a
 line bundle on $Z$, seen as a sheaf on $Y$. The $k$-th order gaussian map
 associated to these data is
\begin{equation}\label{gaussian2}
 \gamma_{L,M_Z}^k:H^0(Y\times
Y,I^k_{\Delta_Y}\otimes L\boxtimes
 M_Z)\rightarrow H^0(Y,S^k{\Omega^1_Y}\otimes L\otimes M_Z).
 \end{equation}
 Sequence (\ref{Ik}), tensored with
 $L\boxtimes M_Z$ remains exact\footnote{To see this, since $q^*M_Z$ is locally free on $Y \times Z$, and $I^k_{\Delta_Y}/I^{k+1}_{\Delta_Y}$ is locally free on $\Delta_Y$, it suffices to show that $tor_1^{\OO_{Y \times Y}}({\mathcal O}_{\Delta_Y}, {\mathcal O}_{Y \times Z})=0.$
One calculates such $tor$ tensoring by $\OO_{\Delta_Y}$ the exact sequence \  \
 $0 \rightarrow I_{Y \times Z} \rightarrow
  {\mathcal O}_{Y \times Y} \rightarrow {\mathcal O}_{Y \times Z} \rightarrow 0,$ \
 obtaining
$$0 \rightarrow    tor_1^{\OO_{Y \times Y}}({\mathcal O}_{\Delta_Y}, {\mathcal O}_{Y \times Z})  \rightarrow  I_{Y \times Z} \otimes  {\mathcal O}_{\Delta_Y} \rightarrow
  {\mathcal O}_{\Delta_Y} \rightarrow {\mathcal O}_{\Delta_Z} \rightarrow 0,$$
  where $\Delta_Z$ is the diagonal of $Z$ in $Y \times Z$. A calculation in local coordinates shows that  $I_{Y \times Z} \otimes  {\mathcal O}_{\Delta_Y}$ is the ideal sheaf of $\Delta_Z$ in $\Delta_Y$. Therefore $tor_1^{\OO_{Y \times Y}}({\mathcal O}_{\Delta_Y}, {\mathcal O}_{Y \times Z})=0$.}.
Hence, as above, the domain of
 these $k$-th gaussian maps is the kernel of the previous ones:
\begin{equation}\gamma^k_{L,M_Z}: ker \gamma^{k-1}_{L,M_Z}\rightarrow
H^0(S^k{\Omega^1_Y}\otimes L\otimes M_Z).
\end{equation} In this
note we will deal with the following setup: $X$  an abelian surface,
$C\subset X$ a smooth and irreducible curve of genus $g\ge 2$. The
line bundles on $X$ and $C$ will be respectively $L= \OO_X(C)$
(necessarily ample) and $K_C$. The curve-surface gaussian maps of
order $\le 2$ associated to these data are the following. The
curve-surface multiplication map
\begin{equation}\label{gamma0}\gamma^0_{X,C}:H^0(X,\OO_X(C))\otimes H^0(C,K_C)\rightarrow
H^0(K_C^2)\end{equation} and the first and second curve-surface
gaussian maps
\begin{equation}\label{gamma1}\gamma^1_{X,C}:H^0(X\times
X,I_{\Delta_X}\otimes \OO_X(C)\boxtimes K_C)\rightarrow
H^0(\Omega^1_X\otimes K_C^2)\end{equation}
\begin{equation}\label{gamma2}\gamma^2_{X,C}:H^0(X\times X,I^2_{\Delta_X}\otimes
\OO_X(C)\boxtimes K_C)\rightarrow H^0(S^2\Omega^1_X\otimes K_C^2)
\end{equation}
\subsection{First gaussian maps and vector bundles}
It is technically useful to see the  gaussian maps
(\ref{gaussian1}) and (\ref{gaussian2}) defined as the $H^0$ of
maps of coherent sheaves on the variety $Y$, rather than on the
cartesian product. This is achieved as follows: let $p$ and $q$
the two projections of $Y\times Y$. Applying ${p}_*$ to the exact
sequences (\ref{Ik}) tensored by $M$ one gets the exact sequences
\begin{equation}\label{vbs}
0\rightarrow {p}_*(I^{k+1}_{\Delta_Y}\otimes q^*M)\rightarrow
{p}_*(I^k_{\Delta_Y}\otimes
q^*M)\buildrel{\varphi^k}\over\rightarrow S^k\Omega^1_Y\otimes M
\end{equation}
 The gaussian maps $\gamma^k_{L,M}$ of (\ref{gaussian1}) are
obtained by tensoring with $L$ and taking $H^0(L\otimes\varphi^k)$.
The same with the gaussian map $\gamma^k_{L,M_Z}$.

Let us spell out how the gaussian maps $\gamma^1_{C}$,
$\gamma^1_{X,C}$, look like in this setting. Let $R_C$ be the the
kernel of the evaluation map of $K_C$:
\begin{equation}
\label{1C} 0 \rightarrow R_C \stackrel{f}\rightarrow H^0(K_C)
\otimes {\mathcal O}_C \rightarrow K_C \rightarrow 0
\end{equation}
(i.e. sequence (\ref{vbs}) for $Y=C$, $M=K_C$, $k=0$). By
(\ref{vbs}) (same setting) for $k=1$ we have the natural map
\begin{equation}\label{g} R_C\buildrel{g}\over\rightarrow
K_C^2. \end{equation}
Tensoring with $K_C$ and taking $H^0$ one obtains the Wahl map
\begin{equation}\label{0C}\gamma^1_C:H^0(R_C\otimes K_C)\rightarrow
H^0(K_C^3).
\end{equation}

Next, let $G$ be the kernel of the evaluation map of $K_C$,
\emph{seen as a sheaf on $X$},
\begin{equation}\label{1X} 0
\rightarrow G \rightarrow H^0(K_C) \otimes {\mathcal O}_X
\rightarrow K_C \rightarrow 0
\end{equation}
(i.e. sequence (\ref{vbs}) for $k=0$, $Y=X$, $M=K_C$). From
(\ref{vbs}) for $k=1$ (same setting) one has the map
\begin{equation}\label{h} G\buildrel{h}\over \rightarrow
\Omega^1_X\otimes K_C.
\end{equation}
It is easily seen that $G$ is a locally free sheaf on $X$, which is
sometimes called a Lazarsfeld's bundle, since this type of
construction was systematically used in \cite{laz}. Tensoring with
$K_C$ and taking $H^0$ one obtains the first curve-surface gaussian
map
\begin{equation}\label{0X}\gamma^1_{X,C}:H^0(G\otimes K_C)\rightarrow
H^0(\Omega^1_X\otimes K_C^2). \end{equation}

\begin{LEM}
Restricting $G$ to $C$, one obtains the exact commutative diagram
$$
\label{curve-surface}
\xymatrix{&&&0\ar[d]&0\ar[d]\\&&&R^2_C\ar[r]^{=}\ar[d]^{\alpha}&R^2_C\ar[d]\\
&0\ar[r] & {\mathcal O}_C\ar[d]^=\ar[r] &  G_{|C}\ar[d]^h\ar[r] &  R_C\ar[r]\ar[d]^g & 0\\
&0 \ar[r] &{\mathcal O}_C  \ar[r]& K_C \otimes \Omega^1_X\ar[r]&
K_C^2 \ar[r]&0}$$ where $R^2_C$ is the kernel of $g$ \emph{(see also
(\ref{R2}) below)}.
\end{LEM}
\proof Restricting sequence (\ref{1X}) to $C$ one has
$$ 0\rightarrow tor_1^{\OO_X}(K_C,\OO_C)
\rightarrow G_{|C} \rightarrow H^0(K_C) \otimes {\mathcal O}_C
\rightarrow K_C \rightarrow 0,$$ where \
$tor_1^{\OO_X}(K_C,\OO_C)\cong K_C\otimes
tor_1^{\OO_X}(\OO_C,\OO_C)\cong \OO_C$.\qed

\subsection{Second gaussian maps and vector bundles}
 Next, the second gaussian maps $\gamma^2_C$ and
 $\gamma^2_{X,C}$.  One has the exact
 sequence
 \begin{equation}
\label{pp} 0 \rightarrow I^{2}_{\Delta_Y} \rightarrow \OO_{Y \times
Y} \rightarrow \OO_{\Delta_Y^2}\rightarrow 0
\end{equation}
where $Y$ is either the surface $X$ or the curve $C$ and
$\Delta_Y^2$ denotes the first infinitesimal neighborhood.

 If $Y=C$, tensoring sequence (\ref{pp}) with $q^*K_C$ and applying
$p_*$  one gets the exact sequence
\begin{equation}\label{fprime}0\rightarrow
R_C^2\buildrel{f^\prime}\over\rightarrow
H^0(K_C)\otimes\OO_C\buildrel{ev}\over\rightarrow P_C(K_C) \ \
 \end{equation}
 where
 \begin{equation}\label{R2}
 R^2_C=p_*(q^*(K_C)\otimes I^{2}_{\Delta_C}),
 \end{equation}
 and
 $$P_C(K_C)= p_*(q^*(K_C)\otimes \OO_{\Delta_C^2})$$ is the
 bundle of principal parts of $K_C$.
 \begin{remark}\label{immersion}The  commutative diagram
$$
\xymatrix{&&0\ar[d]&0\ar[d]\\&&R^2_C\ar[r]^{=}\ar[d]&R^2_C\ar[d]\\
&0\ar[r] & R_C\ar[d]^g\ar[r] & H^0(K_C) \otimes {\mathcal O}_C\ar[d]^{ev}\ar[r] &  K_C\ar[r]\ar[d]^= & 0\\
&0 \ar[r] &K_C^2  \ar[r]& P_C(K_C)\ar[r]&
K_C \ar[r]&0}.$$
shows that the map $g$ of (\ref{g}) is surjective if and only if the evaluation map $ev$ is surjective. This is in turn equivalent to the immersivity of the canonical map, which holds if and only if $C$ is non-hyperelliptic.
\end{remark}

 Sequence (\ref{vbs}) for $k=2$, $Y=C$ and $M=K_C$ provides the natural map
 \begin{equation}\label{gprime}
 R^2_C\buildrel{g^\prime}\over\rightarrow K_C^3.
 \end{equation}
Tensoring with $K_C$ and taking $H^0$ one obtains the second Wahl
map
\begin{equation}\label{2C}\gamma^2_C:H^0(R^2_C\otimes K_C)\rightarrow H^0(K_C^4),
\end{equation}

Finally, let us work out the second curve-surface gaussian map.  Let
us consider the sequence (\ref{pp}) with $Y$ equal to the surface
$X$. Tensoring sequence (\ref{pp}) with $q^*(K_C)$ and applying
$p_*$ one gets
\begin{equation}\label{pponX}
0\rightarrow G^2\rightarrow
H^0(K_C)\otimes\OO_X\buildrel{ev}\over\rightarrow P_X(K_C)
 \end{equation}
 where $G^2=p_*(q^*(K_C)\otimes I^{2}_{\Delta_X})$. The sheaf
  $P_X(K_C)= p_*(q^*(K_C)\otimes \OO_{\Delta_X^2})$ could be referred to as the "the sheaf of
 principal parts of $K_C$ on $X$". Applying $p_*$ to the exact
 sequence
 $$0\rightarrow I_{\Delta_X}/I^2_{\Delta_X}\otimes
 q^*(K_C)\rightarrow \OO_{\Delta^2_X}\otimes
 q^*(K_C)\rightarrow \OO_{\Delta_X}\otimes
 q^*(K_C)\rightarrow 0$$
 one sees that $P_X(K_C)$ sits into the exact sequence (of $\OO_X$-modules)
 $$0\rightarrow \Omega^1_X\otimes K_C\rightarrow
 P_X(K_C)\rightarrow K_C\rightarrow 0.$$
Sequence (\ref{vbs}) with $k=2$, $Y=X$ and $M=K_C$ provides the map
 \begin{equation}\label{hprime}G^2\buildrel{h^\prime}\over\rightarrow
 S^2\Omega^1_X\otimes K_C
 \end{equation}
Tensoring with $K_C$ and taking $H^0$ one obtains the second
curve-surface gaussian map
\begin{equation}\gamma^2_{X,C}:H^0(G^2\otimes K_C)\rightarrow H^0(S^2\Omega^1_X\otimes
K_C^2).
\end{equation}
Let us consider the exact sequence
\begin{equation}\label{sym} 0\rightarrow \Omega^1_X\otimes K_C^{-1}\rightarrow
S^2{\Omega^1_X}_{|C}\rightarrow K_C^2\rightarrow 0
\end{equation}
obtained taking symmetric products in the cotangent sequence. The following Lemma will be useful in the sequel
\begin{LEM} \label{eprime} Assume that $C$ is non-hyperelliptic. Restricting $G^2$ to $C$ one gets the
commutative exact diagram
$$\xymatrix{
0\ar[r] & {\Omega^1_X}_{|C}\ar[d]^=\ar[r] &
G^2_{|C}\ar[d]^{h^\prime}\ar[r] &  R^2_C\ar[r]\ar[d]^{g^\prime}
& 0\\
0 \ar[r] &{\Omega^1_X}_{|C}  \ar[r]& K_C \otimes
S^2{\Omega^1_X}\ar[r]& K_C^3 \ar[r]&0}
$$
\end{LEM}
\proof We have the
  commutative exact diagram
 \begin{equation}\label{big}\xymatrix{&0\ar[d]&0\ar[d]&0\ar[d]\\
0\ar[r]&G(-C)\ar[r]\ar[d]&G^2\ar[r]\ar[d]&R^2_C\ar[r]\ar[d]&0\\
 0\ar[r]&H^0(K_C)\otimes\OO_X(-C)\ar[r]\ar[d]^{ev(-C)}&H^0(K_C)\otimes\OO_X\ar[r]\ar[d]^{ev}
 &H^0(K_C)\otimes \OO_C\ar[r]\ar[d]^{ev}&0\\ 0\ar[r]&\OO_C\ar[r]\ar[d]&P_X(K_C)\ar[r]&P_C(K_C)\ar[r]&0
 \\&0}
 \end{equation}
 Restricting the top row to $C$ one obtains
$$0 \rightarrow tor_1^{\OO_S}(R^2_C,\OO_C)\buildrel{a}\over\rightarrow
G_{|C}\otimes K_C^{-1}\buildrel{b}\over\rightarrow
G^2_{|C}\rightarrow R^2_C\rightarrow 0,$$ where
$tor_1^{\OO_S}(R^2_C,\OO_C)\cong R^2_C\otimes
tor_1^{\OO_S}(\OO_C,\OO_C) \cong R^2_C\otimes K_C^{-1}.$

We claim that the map $a:R^2_C\otimes K_C^{-1}\rightarrow
G_{|C}\otimes K_C^{-1}$ is the map $\alpha$ of the diagram of
Lemma \ref{curve-surface} tensored with $K_C^{-1}$.  Granting this
for the moment, let us conclude the proof. Again in the diagram of
Lemma \ref{curve-surface} we have that, since $C$ is assumed to be
non-hyperelliptic, the map $g$ is surjective (Remark
\ref{immersion}). Therefore also the map $h$ is surjective. By the
Claim, this means that $Im(b)={\Omega^1_X}_{|C}$. This proves the
Lemma. To prove the Claim, we consider the commutative diagram
\begin{equation}\label{big}\xymatrix{
0\ar[r]&R^2_C\otimes K^{-1}_C\ar[r]^{a}&G_{|C}\otimes K^{-1}_C\ar[r]^{b}\ar[d]&G^2_{|C}\ar[r]\ar[d]&R^2_C\ar[r]&0\\
 &&G_{|C}\ar[r]^{=}\ar[d]^{=}&G_{|C}\ar[d]
 &&\\
 0\ar[r]&R^2_C\ar[r]^{\alpha}&G_{|C}\ar[r]^{h\quad}\ar[d]&\Omega^1_{X|C}\otimes K_C\ar[r]\ar[d]&0
 \\&&0&0&&},
 \end{equation}
It follows that $a=\alpha\otimes K_C^{-1}.$ \qed

\begin{remark}\label{hyper} In the same way one proves that, if $C$ is hyperelliptic, restricting $G^2$ to $C$ one
obtains the following commutative exact diagram (notation as in the
above proof)
\begin{equation}
\xymatrix{ 0\ar[r] & Im(h)\otimes K_C^{-1}\ar[d]\ar[r] &
G^2_{|C}\ar[d]^{h^\prime}\ar[r] & R^2_C\ar[r]\ar[d]^{g^\prime}
& 0\\
0 \ar[r] &{\Omega^1_X}_{|C}  \ar[r]& K_C \otimes
S^2{\Omega^1_X}\ar[r]& K_C^3 \ar[r]&0}
\end{equation}
\end{remark}

\section{The first Wahl map of curves on abelian surfaces}
According to Beauville and M\'erindol  (\cite{bm}), we focus on  the extension class
$$e\in {\rm Ext}^1(K_C,K_C^{-1})$$ of the
cotangent sequence \begin{equation}\label{cotext}0\rightarrow
K_C^{-1}\rightarrow {\Omega^1_X}_{|C}\rightarrow K_C\rightarrow
0.\end{equation}
 Since the surface $X$ is abelian, the cotangent
bundle is trivial. This implies  that $e\ne 0$. Via Serre duality,
${\rm Ext}^1(K_C,K_C^{-1})$ is identified to $H^0(K_C^3)^\vee$.
Under this identification we have
\begin{LEM}\label{cotangent} Let $X$ be an abelian surface and let $C\subset X$ be
a smooth and irreducible curve. Assume that the curve-surface
multiplication map
$$\gamma^0_{X,C}:H^0(X,\Oh_X(C))\otimes H^0(C,K_C)\rightarrow
H^0(K_C^2)$$ is surjective. Then $$e\not\in {\rm
Ann}(Im(\gamma^1_C)).$$
\end{LEM}
\proof By functoriality of Serre duality, the Serre-dual  gaussian map
$${\gamma^1_C}^\vee:{\rm
Ext}^1(K_C^2,\Oh_C)\rightarrow {\rm Ext}^1(R_C,\Oh_C)$$ is obtained
applying ${\rm Ext}^1(\cdot,\OO_C)$ to the map $g$ of (\ref{g}). We
will prove that ${\gamma^1_C}^\vee(e)$ is non-zero. To this purpose,
applying ${\rm Ext}^1(\cdot,\OO_C)$ to the map $f$ of sequence
(\ref{1C}), one gets the map
$$\psi: {\rm Hom}(H^0(K_C),H^1(\OO_C))\cong {\rm
Ext}^1(H^0(K_C)\otimes\OO_C,\OO_C)\rightarrow {\rm
Ext}^1(R_C,\OO_C). $$ Now let
$$\delta: H^0(K_C)\rightarrow H^1(\OO_C)$$ be the composition of the
coboundary map of the standard exact sequence
\begin{equation}
\label{22C} 0\rightarrow \OO_X\rightarrow \OO_X(C)\rightarrow
K_C\rightarrow 0
\end{equation}
with the natural injection  $i: H^1(\OO_X) \rightarrow H^1(\OO_C)$.
It follows form its definition that, since the surface $X$ is
abelian, the map $\delta$ is non-zero (its image is $i(H^1(\OO_X))$.
The key point is the following
 \begin{CLA}\label{phipsi}
${\gamma^1_C}^\vee(e)=\psi(\delta).$
 \end{CLA}
 Admitting the Claim for the moment, let us finish the proof.\\
 Applying ${\rm Ext}^1(\cdot,\OO_C)$ to sequence (\ref{1C}), one sees
 that the kernel of $\psi$ is the image of the  map
 \begin{equation}\label{comult} {\gamma^0_{C}}^\vee: {\rm
Ext}^1(K_C,\OO_C)\rightarrow {\rm Hom}(H^0(K_C),H^1(\OO_C)),
\end{equation}
which is the dual of the multiplication map (\ref{symmetric}) for
$X=C$ and $L=K_C$. Therefore, ${\gamma^1_C}^\vee(e)=0$ means, by the
Claim, that $\delta\in\ker(\psi)=Im({\gamma^0_{C}}^\vee)$. The
restriction of the map $\delta$ to the linear series
$$V=Im(H^0(X,\OO_X(C))\rightarrow H^0(K_C))$$
 is zero, by definition.
Hence, if ${\gamma^1_C}^\vee(e)=0$, then there exists an element
$\eta\in{\rm Ext}^1(K_C,\OO_C)$ such that
${\gamma^0}_{C}^\vee(\eta)=\delta$, and $\eta$ belongs to the kernel
of the map
$${\rm Ext}^1(K_C,\OO_C)\rightarrow {\rm Hom}(H^0(\OO_X(C)),H^1(\OO_C)).$$ Now
this last map is the dual of the curve-surface multiplication map
$\gamma^0_{X,C}:H^0(\OO_X(C))\otimes H^0(K_C)\rightarrow
H^0(K_C^2)$. Since $\delta$ is non-zero, $\eta$ is non-zero. But
this in contrast with the assumption that the curve-surface
multiplication map is surjective.

\medskip \emph{Proof of the Claim. } Lemma \ref{curve-surface} shows that
 ${\gamma^1_C}^\vee(e)$ is the extension class of the sequence
 $$0\rightarrow \OO_C\rightarrow G_{|C}\rightarrow R_C\rightarrow 0. $$
Next, we compute $\psi(\delta)$. The exact sequence (\ref{22C})
yields naturally the extension
\begin{equation}
\label{deltasopra} 0\rightarrow \OO_X \rightarrow E \rightarrow
H^0(K_C) \otimes \OO_{X}\rightarrow 0
\end{equation}
sitting in the following commutative and exact diagram

\begin{equation}
\label{E} \xymatrix{
& & &  0\ar[d]& 0\ar[d] & \\
& & &G \ar[r]^=\ar[d] & G \ar[d]& \\
&0\ar[r] & {\mathcal O}_X \ar[d]^=\ar[r] &  E \ar[d]\ar[r] &  H^0(K_C) \otimes \OO_X \ar[r]\ar[d] & 0\\
&0 \ar[r] &{\mathcal O}_X \ar[r]& \OO_{X}(C) \ar[r]\ar[d]& K_C \ar[r] \ar[d]& 0\\
& & & 0&0&}
\end{equation}

This shows that $\delta$ is the extension class of the exact
sequence
\begin{equation}
\label{delta} 0\rightarrow \OO_C \rightarrow E_{|C} \rightarrow
H^0(K_C) \otimes \OO_{C}\rightarrow 0
\end{equation}
obtained by restriction of (\ref{deltasopra}) to $C$. The
commutative diagram with exact rows
\begin{equation}
\label{E_C} \xymatrix{
& 0\ar[r]& \OO_C \ar[r]\ar[d]^=&G_{|C} \ar[r]\ar[d] & R_C\ar[d]^g\ar[r] & 0\\
&0\ar[r] & {\mathcal O}_C \ar[r] &  E_{|C} \ar[r] & H^0(K_C)
\otimes \OO_{C} \ar[r]& 0}
\end{equation}
proves that $\psi(\delta) = (\gamma^1_C)^\vee(e)$, i.e. the Claim. \qed

As a consequence we get the main result of this section (see Theorem
\ref{loose} of the Introduction)
\begin{TEO}\label{firstgaussian}
Let $X$ be an abelian surface and let $C\subset X$ be a smooth and
irreducible curve. Assume that the curve-surface multiplication map
$$\gamma^0_{X,C}:H^0(X,\Oh_X(C))\otimes H^0(C,K_C)\rightarrow
H^0(K_C^2)$$ is surjective.  Then the natural map
$$\coker\gamma_X^1\rightarrow\coker \gamma_C^1$$
is surjective. In particular if, in addition, the Gaussian map $\gamma^1_X$ on the
surface $X$ is surjective then the Gaussian map $\gamma^1_C$ on the
curve $C$ is surjective.
\end{TEO}
\proof We have the natural commutative diagram $$\xymatrix{
\wedge^2H^0(X,\OO_X(C)) \ar[dd]\ar[r]^{\gamma^1_{X} \ \ \ \ } & H^0(\Omega^1_X \otimes {\mathcal O}_X(2C)) \ar[dr]^{q_1} \\
&&H^0(\Omega^1_{X|C} \otimes K_C^2) \ar[dl]^{q_2}\\
\wedge^2H^0(K_C) \ar[r]_{\gamma^1_C}& H^0(K_C^3)}$$
 By  Serre
duality the cotangent extension class $e$ is identified to the
linear functional $f_e:H^0(K_C^3)\rightarrow H^1(K_C)\cong{\mathbb
C}$ defined by the coboundary map of contangent sequence tensored
with $K_C^2$:$$0\rightarrow K_C\rightarrow{\Omega^1_X}_{|C}\otimes
K_C^2\rightarrow K_C^3\rightarrow 0.$$ The map $q_1$ is
surjective. The theorem follows since, by Lemma \ref{cotangent}
the one-codimensional subspace $Im\,(q_2\circ q_1)=Im\, (q_2)$ does
not contain the image of the the Gaussian map
$\gamma^1_C$.\qed

\begin{remark} (a) If $X$ is a K3 surface, instead of  abelian, Claim \ref{phipsi} recovers Beauville-M\'erindol's theorem, asserting the opposite happens, namely that $e\in Ann({\gamma^1_C}^\vee)$, and therefore, (at least if $e\ne 0$) $\gamma^1_C$ is not surjective (\cite{bm}).\\
(b) Let $V$ be the image of the natural map \  $H^0(\OO_X(C))\rightarrow H^0(K_C),$ and let
$$\gamma^1_{V,C}:\wedge^2V\rightarrow H^0(K_C^3)$$ be the restricted Wahl map. By a similar argument, Claim \ref{phipsi} proves that  $e\in Ann(Im \, \gamma^1_{V,C})$. Hence $\gamma^1_{V,C}$ is not surjective.
\end{remark}

 Note that -- by an immediate computation -- the assumption of
Theorem \ref{firstgaussian} holds as soon as one as the
multiplication map on $X$:
$$\gamma^0_X:S^2H^0(X,\OO_X(C))\rightarrow H^0(X,\OO_X(2C))$$ is surjective. This
is known to hold when:

\noindent  -  $\OO_X(C)$ is a power of order at least $3$ of a
(necessarily ample) line bundle (Koizumi, \cite{lb} Th. 7.3.1);

\noindent -    $\OO_X(C)$ is a second power, and no point of the
finite group $K(\OO_X(C))$ is a base point of a symmetric line
bundle algebraically equivalent to $\OO_C(C)$ (Ohbuchi, \cite{lb}
Prop. 7.2.3). This result and the previous one hold for abelian varieties of any
dimension.

\noindent - $\OO_X(C)$ is not a power (i.e. $\OO_X(C)$ is of type
$(1,d)$ ),  it is birational, and $d=g(C)-1\ge 7$, if $d$ is odd,
$d\ge 14$, if $d$ is even, (Lazarsfeld, \cite{laz2}).

\noindent - $\OO_X(C)$ of type $(1,d)$,
the N\'eron-Severi group
 $NS(X)$ is generated by $c_1(\OO_X(C))$ and
 $d\ge 7$ (Iyer, \cite{iy}).

 Finally, let us focus on surjectivity of the first gaussian map on $X$:
 \begin{equation}\label{gaussianX}\gamma^1_X:\wedge^2H^0(X,\OO_X(C))\rightarrow
 H^0(\Omega^1_X(2C)).
 \end{equation}
 Unlike for multiplication maps, effective surjectivity criteria for gaussian maps of polarizations of type $(1,d)$ on abelian
 surfaces are not available at present. The only  result we are aware of is the analogue of Koizumi's theorem,
 asserting that, if $\OO_X(C)$ is at least a $5$-th power, then the
 gaussian map (\ref{gaussianX}) is surjective (\cite{par} Th.2.1).
 Therefore, as a consequence of Theorem \ref{firstgaussian}, we have
 \begin{COR}\label{surj} Let $X$ be an abelian surface and  let ${\mathcal L}$
 be an ample line bundle on $X$, and let $k\ge 5$. For all smooth
 and irreducible curves $C\in |{\mathcal L}^k|$ the Wahl map
 $$\gamma^1_C:\wedge^2H^0(K_C)\rightarrow H^0(K_C^3)$$
 is surjective.
 \end{COR}
 It is worth to note that suitable surjectivity criteria for gaussian maps
 on abelian surfaces -- analogous to those due to Lazarsfeld and
 Iyer for multiplication maps -- would imply, as in the previous Corollary, the surjectivity
 of the Wahl map of general
 curves of any suitably high genus lying on abelian surfaces, thus
 providing a "without degeneration" proof of the theorem of
 Ciliberto-Harris-Miranda.

\section{The second Wahl map of curves on abelian surfaces}
The present Section is entirely devoted to the proof of Theorem A.

Tensoring the symmetric square of the cotangent sequence (\ref{sym})
with $K_C^2$ one gets
\begin{equation} 0\rightarrow \Omega^1_X\otimes K_C\cong K_C^{\oplus 2}\rightarrow
S^2{\Omega^1_X}_{|C}\otimes K_C^2\rightarrow K_C^4\rightarrow 0
\end{equation}
whose coboundary map  $$f_{e^\prime}: H^0(K_C^4)\rightarrow
H^1(\Omega^1_X\otimes K_C)\cong H^1(K_C)^{\oplus 2}\cong{\mathbb
C}^{\oplus 2}$$ is identified, by Serre duality, to the extension
class $e^\prime\in{\rm Ext}^1(K_C^3,{\Omega^1_X}_{|C})$ of
 sequence
(\ref{sym}). The first part of the statement of Theorem \ref{main}
is equivalent to:
\begin{equation}\label{stat1}
f_{e^\prime}\circ \gamma^2_C= 0. \end{equation}
 As for the case
of first gaussian maps, it is easier to work in the dual setting.
Applying ${\rm Ext}^1(\cdot, {\Omega^1_X}_{|C})$ to the map
$g^\prime$ of (\ref{gprime}) one gets the map
$$\phi:{\rm Ext}^1(K_C^3,{\Omega^1_X}_{|C})\rightarrow {\rm
Ext}^1(R^2_C,{\Omega^1_X}_{|C})$$ (which is identified to two copies
of the dual map of the second Wahl map) and it is easily seen that
(\ref{stat1}) is equivalent to the fact that
\begin{equation}\label{stat2}\phi(e^\prime)=0.
\end{equation}
Applying ${\rm Ext}^1(\cdot,\Omega^1_{X|C})$ to the map $f^\prime$
of (\ref{fprime}) we get the map
$$\psi: {\rm Hom}(H^0(K_C),H^1(\Omega^1_{X|C}))={\rm
Ext}^1(H^0(K_C)\otimes\OO_C,\Omega^1_{X|C})\rightarrow {\rm
Ext}^1(R^2_C,\Omega^1_{X|C}).$$
 Now let us denote by $\tilde{\delta}$ the composition of the
coboundary map $H^0(K_C)\rightarrow H^1(\OO_X)$ of the standard
exact sequence (\ref{22C}), and the  map $H^1(d): H^1(\OO_X)
\rightarrow H^1(\Omega^1_{X|C})$, induced by the derivation $d:
\OO_X \rightarrow \Omega^1_{X}$. Note that $H^1(d)$ is the zero
map since the Hodge-Fr\"olicher spectral sequence degenerates at level of $E_1$.

\begin{CLA}
\label{phipsi2} $\phi(e^\prime)=\psi(\tilde{\delta}) =0$
\end{CLA}
The first part of Theorem \ref{main}, i.e. (\ref{stat1}), follows
immediately since $\tilde\delta=0$.

\medskip\emph{Proof of the Claim. } Assume that $C$ is non-hyperelliptic. Lemma \ref{eprime} shows
that $\phi(e^\prime)$ is the class of the sequence $$0\rightarrow
{\Omega^1_X}_{|C}\rightarrow G^2_{|C}\rightarrow R^2_C\rightarrow
0$$ Next, we compute $\psi(\tilde\delta)$. The zero map
$H^1(d) \in {\rm Hom}(H^1(\OO_X),H^1(\Omega^1_{X}))$ can be seen as
the class of
 the extension
\begin{equation}
\label{H} 0 \rightarrow \Omega^1_{X} \rightarrow H:=
R^1q_*(I_{\Delta_X}^2) \rightarrow  H^1(\OO_X)
 \otimes \OO_X \rightarrow 0
\end{equation}
obtained applying $q_*$ to the exact sequence (\ref{Ik}) with $k
=1$. In fact, by Leray spectral sequence the coboundary map of (\ref{H})
can be identified with the map in cohomology:
$$H^1(X,{\mathcal O}_X)\stackrel{^F_{\sim}}\rightarrow H^1(X\times X, I_{\Delta_X})\stackrel{G}\rightarrow H^1(X\times X,
I_{\Delta_X}/I^2_{\Delta_X})\cong H^1(\Omega^1_X)$$ induced by the
exact sequence (\ref{Ik}) with $k =1$.
K\"unneth formula gives the isomorphism $H^1(X,{\mathcal O}_X) \oplus H^1(X,{\mathcal O}_X) \stackrel{\sim} \rightarrow H^1(X\times X, \OO_{X \times X}),$ $(\alpha,\beta) \mapsto p^*(\alpha) + q^*(\beta)$.
Therefore the map $F: H^1(X,{\mathcal O}_X) \rightarrow H^1(X\times X, I_{\Delta_X})$ is given by $\alpha \mapsto p^*(\alpha) - q^*(\alpha)$, hence $G \circ F = H^1(d)$, since the map $d: \OO_X \rightarrow I_{\Delta_X}/I^2_{\Delta_X}$ is exactly the map $\alpha \mapsto p^*(\alpha) - q^*(\alpha)$.

We have the natural exact
diagram
\begin{equation}
\label{e} \xymatrix{
&& 0\ar[r]& R^2_X\ar[r]\ar[d]&H^0(\OO_X(C))\otimes \OO_X \ar[r]\ar[d] & P_X(\OO_X(C)) \ar[d]& \\
&&0\ar[r]& G^2\ar[d]\ar[r] &  H^0(K_C)\otimes \OO_X \ar[d]\ar[r] &  P_X(K_C)\ar[d] & \\
&0\ar[r]&\Omega^1_X\ar[r] &H\ar[r]& H^1(\OO_X)\otimes \OO_X\ar[r]&
0 &}
\end{equation}
where $R^2_X=q_*(I^{2}_{\Delta_X}\otimes p^*({\mathcal O}_X(C)))$.
The first column of the diagram is obtained from the exact
sequence
$$0 \rightarrow p^*(\OO_X)  \rightarrow p^*(\OO_X(C))     \rightarrow  p^*(K_C) \rightarrow 0$$
tensoring with $I^2_{\Delta_X}$ and applying $q_*$. Notice that
tensoring with $I^2_{\Delta_X}$, the above exact sequence remains
exact, since
$$tor_1^{\OO_{X \times X}}(p^*(K_C), I^2_{\Delta_X}) = tor_2^{\OO_{X \times X}}(p^*(K_C), I_{\Delta_X}/I^2_{\Delta_X})=0.$$

Restricting  the two bottom rows to $C$ one gets
\begin{equation}
 \xymatrix{
&0\ar[r]& \Omega^1_{X|C}\ar[r]\ar[d]^{=}& G^2_{|C}\ar[d]\ar[r] &  R^2_C \ar[d]^{\lambda}\ar[r] &  0 & \\
&0\ar[r]&\Omega^1_{X|C}\ar[r] &H_{|C}\ar[r]& H^1(\OO_X)\otimes
\OO_C\ar[r]& 0 &}
\end{equation}
where $\lambda$ is the composition
$$R^2_C \stackrel{g} \rightarrow  H^0(K_C)\otimes \OO_C\stackrel{\delta \otimes \OO_C} \rightarrow H^1(\OO_X)\otimes \OO_C.$$
Hence $\phi(e^\prime)=\psi(\tilde\delta)$. This proves Claim
\ref{phipsi2} if $C$ is non-hyperelliptic. If $C$ is hyperelliptic
the same argument applies, using the diagram of Remark \ref{hyper}
instead of the one of Lemma \ref{eprime}. This concludes the  proof
of Claim \ref{phipsi2}, and hence of the first part of the statement
of Theorem \ref{main}.

The last part of the statement follows as in Theorem
\ref{firstgaussian} from the commutative diagram
\begin{equation}\xymatrix{
I_2(\OO_X(C))\ar[dd]\ar[r]^{\gamma^2_X \ \ \ \ \ \ } & H^0(S^2 \Omega^1_X \otimes {\mathcal O}_X(2C)) \ar[dr] \\
&&H^0(S^2 \Omega^1_{X|C} \otimes K_C^2) \ar[dl]\\
I_2(K_C) \ar[r]_{\gamma^2_C}& H^0(K_C^4)}
\end{equation}
\qed

 Finally, it is known  that, if $\OO_X(C)$ is at least a $7$-power of
a (necessarily ample) line bundle on $X$, then the second gaussian
map $\gamma^2_{X}$ is surjective (this
 is part of a general
result on surjectivity of higher gaussian maps of any order for
powers of ample line bundles on abelian varieties, \cite{par} Th.
2.2). Hence we have the following
\begin{COR} Let $X$ be an abelian surface, let ${\mathcal L}$ be an ample line bundle on $X$ and let $k\ge 7$. Then, for every smooth and
irreducible curve $C\in|{\mathcal L}^k|$, the image of second Wahl
map
$$\gamma^2_C:I_2(K_C)\rightarrow H^0(K_C^4)$$
is the 2-codimensional subspace $S^2H^0(\Omega^1_X)\cdot
H^0(K_C^2)$.
\end{COR}

\end{document}